\documentclass[reqno]{amsart}
\usepackage{epsfig,graphicx,amssymb,amscd}
\DeclareMathOperator{\dm}{dom}
\DeclareMathOperator{\im}{img}
\DeclareMathOperator{\spc}{spec}

\newcommand{\rl}{{\mathbb{R}}}
\newcommand{\cx}{{\mathbb{C}}}
\newcommand{\dbar}{\overline{\partial}}

\newcommand{\bdry}{\mathsf{b}}

\newcommand{\csor}{{\widehat{\otimes}}}
\newcommand{\tensor}{\otimes}

\newcommand{\abs}[1]{\left|{#1}\right|}
\newcommand{\norm}[1]{\left\|{#1}\right\|}

\newtheorem{thm}{Theorem}[section]
\newtheorem{lem}[thm]{Lemma}
\newtheorem{prop}[thm]{Proposition}
\newtheorem{cor}[thm]{Corollary}

\newcommand{\h}{\mathsf{H}}

\begin{document}
\title{Spectrum of the complex Laplacian on  product domains}
 \author{Debraj Chakrabarti}
 \address{Department of Mathematics,  University of Notre Dame,
 Notre Dame, IN 46556, USA} \email{dchakrab@nd.edu}
 \begin{abstract} 
 We show that the spectrum of the complex Laplacian $\Box$ on a product
 of hermitian manifolds is the Minkowski sum of the spectra of the complex Laplacians
 on the factors. We use this to show that the range of the Cauchy-Riemann operator $\dbar$
is closed on a product manifold, provided it is closed on each factor manifold.
 \end{abstract}
\maketitle
\section{Introduction}
The study of the $\dbar$ and $\dbar$-Neumann problems  on product domains
raises a series of interesting questions, which have been studied by
many authors \cite{eh1,eh2,eh3,jb,chsh,fu}.
In \cite{fu}, the method of separation of variables was used to compute the 
spectrum of the complex Laplacian  $\Box=\dbar\dbar^*+\dbar^*\dbar$ on a polydisc,
and for each eigenvalue, the corresponding eigenspace was identified. In this note,
inspired by \cite{fu},
we use the same  technique to compute the spectrum of $\Box$ on an arbitrary product manifold
in terms of  the spectra of $\Box$ on the factors,
and give a description of the spectral representation of $\Box$.
 As a consequence, we obtain a strengthened version (see Theorem~\ref{thm-zucker} below)
of the closed-range theorem  for the $\dbar$ operator on products.

In this note an {\em operator} $A$ on a Hilbert space $\h$ 
(always  assumed  separable)
is a linear map with target $\h$ and source a linear subspace $\dm(A)$ of $\h$. 
The $\Box$-Laplacian on a hermitian manifold $\Omega$ is a densely defined
selfadjoint linear operator on the Hilbert space $L^2_*(\Omega)$ of
$L^2$ differential forms on $\Omega$. Further, $\Box$ maps the subspace $L^2_{p,q}(\Omega)$ 
of forms of bidegree $(p,q)$ to itself, and we denote the restriction of $\Box$ to $L^2_{p,q}(\Omega)$
by $\Box_{p,q}$. More details on $\Box$ may be found in $\S$\ref{sec-box} below.

Recall that  the {\em spectrum}  of  $A$ consists of
those $z\in\cx$ such that the operator $A-zI$ does not have a bounded inverse.
We denote by $\sigma(\Omega)$ (resp. $\sigma_{p,q}(\Omega)$) the spectrum of the operator
$\Box$ on $L^2_*(\Omega)$ (resp. $\Box_{p,q}$ on $L^2_{p,q}(\Omega)$).  Since $\Box$ is a 
nonnegative selfadjoint operator it follows each of $\sigma(\Omega)$ and $\sigma_{p,q}(\Omega)$ 
is a closed nonempty subset of the set of  nonnegative real numbers. For non-compact $\Omega$,
the set $\sigma(\Omega)\subset \rl$
need not be discrete: for example, Lemma~\ref{lem-tfae} below implies that for the ``bumped shell" domain
described on pp. 75-76 of \cite{fk},   $0\in\rl$ is a point of accumulation of the spectrum.

For sets of numbers $P_1,\dots, P_N$ we  denote by  $P_1+\dots+P_N$ their {\em Minkowski sum},
i.e. the set  $\{p_1+\dots+p_N\mid p_1\in P_1,\dots, p_N\in P_N\}$. We have the following:
\begin{thm}\label{thm-specsum} For $j=1,\dots, N$, let $\Omega_j$ be a hermitian manifold, and let
$\Omega= \Omega_1\times\dots\times\Omega_N$ be the product manifold with the product hermitian 
metric. Then \begin{equation}\label{eq-sigma}
\sigma(\Omega) = \sigma(\Omega_1)+\dots+\sigma(\Omega_N).
\end{equation}

Further, we have 
\begin{equation}\label{eq-sigmapq} \sigma_{{\boldsymbol{p}},{\boldsymbol{q}}}(\Omega) =\bigcup_{\substack{\sum_{j=1}^N p_j={\boldsymbol p}\\\sum_{j=1}^N q_j={\boldsymbol q}}}
\left(\sigma_{p_1,q_1}(\Omega_1)+\dots+\sigma_{p_N,q_N}(\Omega_N)\right).\end{equation}
\end{thm}
The computation of the spectrum of $\Omega$ described above can be used to deduce 
properties of $\Omega$ from those of $\Omega_1,\dots,\Omega_N$. We consider one example of this approach.
A very important property that the  operator  $\dbar:L^2_*(\Omega)\rightarrow L^2_*(\Omega)$ might possess 
is that of having closed range, i.e., $\im(\dbar)\subset L^2_*(\Omega)$ is a closed subspace. This is equivalent
to the solvability of the $\dbar$-equation $\dbar u=f $ (where $\dbar f=0$, and $f$ is orthogonal to the harmonic forms)
in the $L^2$-sense (see \cite{cs} for details.) The closed range
property holds under suitable convexity assumptions, e.g. when $\Omega$ is pseudoconvex. We can  define the $L^2$-Dolbeault cohomology space
 \[ {H}^*_{L^2}(\Omega) = \frac{\ker(\dbar)}{\im(\dbar)}.\]
 When the closed-range
property holds, this is a Hilbert space with the quotient norm. This space represents 
the obstruction to solving the $\dbar$-problem in the $L^2$-sense on $\Omega$. We have the following:
\begin{thm}\label{thm-zucker} Let $\Omega_1,\dots,\Omega_N$ and $\Omega$ be as in Theorem~\ref{thm-specsum}.
Suppose that for each $j$, the $\dbar$-operator  on $\Omega_j$ in the $L^2$-sense has closed range in $L^2_*(\Omega_j)$.
Then $\dbar: L^2_*(\Omega)\rightarrow L^2_*(\Omega)$ also has closed range.
Furthermore,  the  K\"{u}nneth formula holds for the $L^2$ cohomology:
\begin{equation}\label{eq-kun}   H^*_{L^2}(\Omega)= H^*_{L^2}(\Omega_1)\csor \dots\ \csor H^*_{L^2}(\Omega_N),\end{equation} 
where   $\csor$ denotes the Hilbert space tensor product (cf. $\S$\ref{sec-tensor} below.) 
\end{thm}

The analog of \eqref{eq-kun} for the $L^2$ de Rham cohomology, in the special case  when the cohomology
spaces are finite dimensional,  goes back to the work of Cheeger (see  \cite{chgr2} and especially \cite[p.~614]{chgr1}.)
Another approach, using an explicit solution of the $d$-equation on a product domain was given
by Zucker in  \cite{zuck}. This can be extended to the $\dbar$-equation (see \cite{chsh}),
or, in another direction,  to  abstractly defined products of  ``Hilbert Complexes" (see \cite{brun}.)
A crucial feature in this approach is the assumption that the $d$ or $\dbar$ operator
on the product satisfies a ``Leibniz rule" in the strict operator sense (cf. assumption (i) in the statement
of   \cite[Theorem~2.29]{zuck}.) In practice, this means that some boundary-smoothness or completeness assumptions
must made on the factors $\Omega_j$ in order to get a handle on the domain of the $d$- or $\dbar$- operator (via
Friedrichs' lemma when the $\Omega_j$  have boundaries.) A version of Theorem~\ref{thm-zucker} is proved in \cite{chsh},
where it is assumed that each of the factors
$\Omega_j$ has Lipschitz boundary. Consequently, using the results of \cite{chsh}, we cannot conclude, for example,
that $\dbar$ has closed range on the product domain in $\cx^4$ given as
\[ \Omega= \{z\in\cx^2\mid 0<\abs{z_1}<\abs{z_2}<1\}\times \{z\in\cx^2\mid 1<\abs{z}<2\},\]
although, in both factors $\dbar$ has closed range : the first (the ``Hartogs triangle") is pseudoconvex (but has a singular non-Lipschitz boundary), and for the second we can see \cite{mcs1, mcs2}. Our result above shows that
$\dbar$ has closed range  on $\Omega$, since Theorem~\ref{thm-zucker} 
involves no assumptions on the factor manifolds $\Omega_j$ except
the closed-range property for $\dbar$.

{\em Acknowledgement:} The author wishes to express his gratitude to Professor Mei-Chi Shaw
for  her support and advice, and to the referee for helpful suggestions.

\section{Some results from Functional Analysis}
We recount here some facts from functional analysis which will be used in the proofs
of Theorems~\ref{thm-specsum} and \ref{thm-zucker}. Most of what we need can be found in
e.g. \cite{rs}, but we discuss the required results to set up notation and for completeness. 

\subsection{Multiplication Operators}\label{sec-multop}
Consider the measure space $(X,\mu)=(X,\mathcal{S},\mu)$ (we systematically suppress the $\sigma$-algebra
from the notation from now on), and a real valued measurable function $h$ on $X$ which is finite a.e. (with respect 
to  $\mu$; this is the last time we will mention the measure with ``a.e.") Define the multiplication operator $T_h$ to be the operator on $L^2(X,\mu)$  on the domain
\[ \dm(T_h) = \{f\in L^2(X,\mu)\mid hf \in L^2(X,\mu)\} \]
defined by 
$ T_h f = hf.$ 
Note that if $h=\tilde{h}$ a.e., then the operators $T_h$ and $T_{\tilde{h}}$ on
$L^2(X,\mu)$ are identical.  If $\mu$ is a finite measure, $T_h$ is densely defined. In fact, it is not difficult to see
(cf. \cite[Prop.2, p.260]{rs}) that if the function $h\in L^p(X,\mu)$, then any dense linear subspace of $L^q(X,\mu)$
is a core of $T_h$, provided $p^{-1}+q^{-1}= \frac{1}{2}$. (Recall that for an operator $A$ on a Hilbert space 
$\h$, a {\em core} of $A$  is a linear subspace of $\dm(A)$ which is dense in the graph norm $x\mapsto \norm{x}_{\h} +\norm{Ax}_{\h}$ in $\dm(A)$.)

It is not difficult to see that the operator $T_h$ is selfadjoint, and the spectrum $\spc (T_h)$ of
the operator $T_h$ is identical to the {\em essential range} of the  function $h$.
Recall that the {\em essential range} of a real-valued function $h$  on $(X,\mu)$ is the set
of $\lambda\in\rl$ such that for all $\epsilon>0$, we have
\[ \mu\left\{x\in X\mid \lambda-\epsilon < h(x)< \lambda+\epsilon\right\}>0.\]
Clearly, the essential range is a closed subset of the real line.
Further, $\lambda$ is an eigenvalue of $T_h$ if and only if $\mu(h^{-1}(\lambda))>0$. The corresponding
eigenspace of $T_h$ is the closed subspace of $L^2(X,\mu)$ consisting of functions which vanish a.e.
outside the set  $h^{-1}(\lambda)\subset X$. Considering the special case when $\lambda=0$,
we obtain a natural identification
\begin{equation}\label{eq-eigen}
\ker(T_h)\cong L^2(h^{-1}(0),\mu)\hookrightarrow L^2(X,\mu),
\end{equation}
where the measure space $(h^{-1}(0),\mu)$ is defined by restriction,
and the inclusion of  $L^2(h^{-1}(0),\mu)$ in $L^2(X,\mu)$ is induced by
extension of functions by 0 from $h^{-1}(0)$ to $X$.
When $h^{-1}(0)$ is the empty set, we will define $L^2(h^{-1}(0),\mu)$ to be the 
trivial vector space $\{0\}$. With this understanding \eqref{eq-eigen} is correct
for all multiplication operators $T_h$. 

\subsection{The spectral theorem} \label{sec-spthm}
The spectral theorem, a structure theorem for selfadjoint operators on a Hilbert space,
can be stated in various equivalent forms
(see  \cite[Theorems VIII.4, VIII.5 and VIII.6. pp. 260-264]{rs}; for bounded operators, see the masterly exposition
\cite{hal}) We will use it in the following form: {\em the multiplication operators $T_h$ defined in $\S$\ref{sec-multop} are }(up to an isometric identification of Hilbert spaces) {\em the {\bf only} examples of selfadjoint operators.} More precisely, let
$A$ be a selfadjoint operator on a Hilbert space $\h$ with domain $\dm(A)$.
Then there is a measure space $(X,\mu)$ with $\mu$ a {\em finite} measure, a unitary 
operator (i.e. isometry of Hilbert spaces) $U:\h\rightarrow L^2(X,\mu)$ and a real valued
$h$ on $X$ finite a.e., so that $\dm(A)= U^{-1}(\dm(T_h))$, and
$A= U^{-1} T_h U$. Note that there is no uniqueness here for the space $X$, the measure $\mu$ or the multiplying
function $h$. We will refer to $T_h$ as a representation of $A$ by a multiplication.

We note  that we can without loss of generality assume that the function $h$ in the  conclusion 
of the spectral theorem belongs to $L^p(X,\mu)$ for every $p\geq 1$. Indeed, if this is not already 
the case, we replace the measure space $(X,\mu)$ by a new measure space $(Y,\nu)$, where $Y=X$, and 
$d\nu=e^{-h^2}d\mu$. Note that $f\mapsto e^{\frac{h^2}{2}}f$ defines an isometry from $L^2(X,\mu)$ to $L^2(Y,\nu)$.
On $Y=X$, the operator $A$ is still represented by multiplication by the same function $h$, but now $h\in L^p(Y,\nu)$ for
each $p\geq 1$.

\subsection{Tensor Products}\label{sec-tensor}
Let $\mathsf{H}_1$ and $\mathsf{H}_2$ be  complex vector spaces.  We denote by
$\mathsf{H}_1\tensor \mathsf{H}_2$  the {\em algebraic} tensor
product (over $\cx$) of  $\mathsf{H}_1$ and $\mathsf{H}_2$ : then $\mathsf{H}_1\tensor \mathsf{H}_2$ can be thought of as
the space of finite sums of elements of the type $x\tensor y$,
where $x\in \mathsf{H}_1$ and $y\in \mathsf{H}_2$,  where
 $\tensor: \mathsf{H}_1\times \mathsf{H}_2\rightarrow \mathsf{H}_1\tensor \mathsf{H}_2$ is the
 canonical bilinear map.
If $\h_1$, $\h_2$ are vector spaces of functions defined on spaces
$X_1$, $X_2$ respectively, then the algebraic tensor product $\h_1\tensor \h_2$ can be 
concretely realized as a space of functions on $X_1\times X_2$ by the correspondence
\begin{equation}\label{eq-conten} (f\tensor g)(x_1,x_2)= f(x_1)g(x_2),\end{equation}
followed by linear extension. We will always make this identification.

When $\h_1$ and $\h_2$ are Hilbert spaces. 
We can define an inner product on $\mathsf{H}_1\tensor \mathsf{H}_2$
by setting
\[ (x\tensor y, z\tensor w)= (x,  z)_{\mathsf{H}_1}(y, w)_{\mathsf{H}_2},\]
and extending bilinearly. This is well-defined thanks to the bilinearity of $\tensor$.
This makes $\mathsf{H}_1\tensor \mathsf{H}_2$ into a pre-Hilbert space,
and its completion is a Hilbert space denoted by $\mathsf{H}_1\csor \mathsf{H}_2$, the {\em Hilbert tensor product}
of the spaces $\mathsf{H}_1$ and $\mathsf{H}_2$. The algebraic tensor product $\mathsf{H}_1\tensor \mathsf{H}_2$ sits
inside $\mathsf{H}_1\csor \mathsf{H}_2$ as a dense subspace. 

For $j=1,2$, let $\h_j = L^2(X_j,\mu)$, and let  $\mu_1\tensor \mu_2$ denote the product 
measure on $X_1\times X_2$. Then the injective map
$L^2(X_1,\mu_1)\tensor L^2(X_2,\mu_2)\hookrightarrow L^2(X_1\times X_2, \mu_1\tensor\mu_2)$ 
given by \eqref{eq-conten} extends to an isometric isomorphism  of Hilbert spaces
\begin{equation}\label{eq-natten} L^2(X_1,\mu_1)\csor L^2(X_2,\mu_2) \cong L^2(X_1\times X_2, \mu_1\tensor \mu_2),\end{equation}
and we will always make this identification.

If $S_1$ and $S_2$ are operators on $\h_1$ and $\h_2$ with dense domains $\dm(S_1)$ and $\dm(S_2)$, we can 
define an operator $S_1\tensor S_2$ on $\h_1\csor \h_2$ with dense domain the algebraic tensor product $\dm(S_1)\tensor \dm(S_2)$ 
by setting on the simple tensors $x\tensor y$:
\begin{equation}\label{eq-tenproddef} (S_1\tensor S_2)(x\tensor y) = S_1 x\tensor S_2 y,\end{equation}
and extending bilinearly. 

\subsection{The operator $A_1\tensor I_2+I_1\tensor A_2$}\label{sec-AB}
Let $\h_1$ and $\h_2$ be separable Hilbert spaces, $A_1,A_2$ be densely defined selfadjoint
operators on $\h_1,\h_2$ respectively, 
and let  $I_1, I_2$ respectively denote the identity maps on $\h_1,\h_2$.
We recall here (see \cite[Theorem VIII.3]{rs} or \cite[Theorem~4.14]{tes})
the spectral representation of the operator $B$ on $\h_1\csor \h_2$,  given by
\[ B= A_1\tensor I_2+ I_1\tensor A_2,\]
which is densely defined with domain $\dm(A_1)\tensor\dm(A_2)$.
 For $j=1,2$ let $(X_j,\mu_j)$ be 
measure spaces and  $U_j:\h_j\rightarrow L^2(X_j,\mu_j)$
be unitary isomorphisms of Hilbert spaces given by the spectral theorem of $\S$\ref{sec-spthm},
such that $A_j=U_j^{-1} T_{h_j}U_j$, where $h_j$ is a real valued function on $X_j$,
such that for each $p\geq 1$, we have $h_j\in L^p(X_j,\mu_j)$. Let 
$U=U_1\csor U_2$, so that $U$ is an unitary isomorphism from 
$\h_1\csor \h_2$ onto $L^2(X_1\times X_2,\mu_1\tensor \mu_2)$. Consider the operator 
$A$ on $\h_1\tensor \h_2$ defined by
\begin{equation}\label{eq-def-A} A= U^{-1} T_h U,\end{equation}
where $h$ is the function on $X_1\times X_2$ defined by $h(x_1,x_2)=h_1(x_1)+h_2(x_2)$. 
Observe that, thanks to the hypotheses on $h_1$ and $h_2$, we have for each $p\geq 1$, that $h\in L^p(X_1\times X_2,\mu_1\tensor\mu_2)$.
Then 
by the results recalled in $\S$\ref{sec-multop}, $A$ is a selfadjoint operator on $\h_1\csor \h_2$. 

\begin{lem}
The restriction of $A$ to $\dm(A_1)\tensor \dm(A_2)$ coincides with $B$.
\end{lem}
\begin{proof}
Note that by linearity, we only need to check this on  simple tensor products of the type $f\tensor g$. The proof
is completed by a direct computation, using the fact that $U^{-1}= U_1^{-1}\tensor U_2^{-1}$ on the algebraic 
tensor product $L^2(X_1,\mu_1)\tensor L^2(X_2,\mu_2)$.
\end{proof}

Several important consequences follow from the lemma above:
\newcounter{property}
\begin{list}{\arabic{property}$^\circ$}{\leftmargin=0cm \itemindent=1cm}
\usecounter{property}
\item {\em The operator $B$ is essentially selfadjoint.}

Recall that an {\em essentially selfadjoint} operator is one whose closure is selfadjoint. 
It is easy to see that such an operator has a unique selfadjoint extension, namely, its closure.
It follows that {\em the operator $A$ is the closure of
the operator $B$.}

Since  $A$ is a selfadjoint extension of 
$B$, to show that $B$ is essentially selfadjoint we need to show that
$\dm(B)=\dm(A_1)\tensor\dm(A_2)$ is a core of $A$. Using \eqref{eq-def-A},
and translating the problem to the representation by multiplication operators,
we need to show that  $\dm(T_{h_1})\tensor \dm(T_{h_2})$ is a core of $T_h$.
Since the function $h\in L^4(X_1\times X_2,\mu_1\tensor\mu_2)$ ( in fact $h\in L^p$ for each $p\geq 1$),
it follows (cf. \cite[Prop.~2, p.~260]{rs}) that any  dense linear subspace of
$L^4(X_1\times X_2,\mu_1\tensor\mu_2)$  is a core of $T_h$.
Therefore, to prove the result it is sufficient to show that $\dm(T_{h_1})\tensor \dm(T_{h_2})$ is dense
in $L^4(X_1\times X_2,\mu_1\tensor\mu_2)$.
Since $h_j\in L^p(X_j,\mu_j)$, for all $p\geq 1$ it follows that all simple functions 
(linear combinations of characteristic 
functions of measurable sets) are in $\dm(T_{h_j})$. It follows that the linear  span $\mathfrak{S}$ of characteristic functions of rectangles with measurable sets as edges is  contained in $\dm(T_{h_1})\tensor \dm(T_{h_2})$. But it is
well-known that $\mathfrak{S}$ is dense in $L^p(X_1\times X_2,\mu_1\tensor \mu_2)$ for each $p>0$.

\item The same method of proof can be used to prove a slightly stronger statement: 

{\em If  for $j=1,2$, the linear space
$D_j\subset \h_j$ is a core of $A_j$, then $D_1\tensor D_2$ is a core of $A$.} For details see \cite[Theorem VIII.3]{rs}.

\item Denote by ${\rm{ess. ran}}(f)$ the essential range of a function $f$ (cf. $\S$\ref{sec-spthm}.) Then
for our functions $h_1, h_2, h$, it is easily verified that
\[ {\rm{ess. ran}}(h) =\overline{{\rm{ess. ran}}(h_1)+{\rm{ess. ran}}(h_2)},\]
where the bar denotes closure in the topology of $\rl$. It follows that, the
spectra of $A_1, A_2$ and $A$ are related by
\[ \spc(A) = \overline{\spc(A_1) +\spc(A_2)}.\]
We note here that the set $\spc(A_1)+\spc(A_2)$ need not be closed. Indeed, 
it is easy to construct selfadjoint  operators  $A_1$ and $A_2$  ( both e.g. on the space $\ell^2$ of 
square-summable sequences) 
such that $ \spc(A_1)=\mathbb{N}_+$, the set of positive integers, and $\spc(A_2)=\left\{-\nu-\frac{1}{\nu}\mid \nu\in \mathbb{N}_{+} \right\}$.

\end{list}

\section{The Complex Laplacian}
\label{sec-box}
\subsection{Definition and Basic Properties}
We now recall the definition and basic properties of the complex Laplacian $\Box$. Details
may be found in the texts \cite{fk,cs}. 

Let $\Omega$ be a hermitian manifold, i.e., a complex manifold with a hermitian metric.
We let $L^2_{p,q}(\Omega)$ be the Hilbert space of square integrable differential forms 
of bidegree $(p,q)$ and let $L^2_*(\Omega)$ be the orthogonal Hilbert space sum of the $L^2_{p,q}(\Omega)$,
so that $L^2_*(\Omega)$ is the Hilbert space of all square integrable forms on $\Omega$.

We can define a realization of the $\dbar$-operator as a densely defined closed Hilbert space
operator from the space $L^2_*(\Omega)$ to itself.
This realization has the domain 
\[ \dm(\dbar)= \{f\in L^2_*(\Omega)\mid \dbar f\in L^2_*(\Omega)\}.\]
where $\dbar f$ is taken in the sense of distributions. We denote by $\dbar^*$ the Hilbert space
adjoint of $\dbar$. This is again a densely defined closed operator on $L^2_*(\Omega)$, whose domain  
$\dm(\dbar^*)$ is very different from that of $\dbar$. We define the {\em complex Laplacian} on $\Omega$
to be the operator
\[ \Box = \dbar\dbar^*+\dbar^*\dbar.\]
Then, it can be shown that $\Box$ is a densely defined closed and unbounded operator,
which is {\em selfadjoint} and {\em nonnegative.}  Note that by the definition of domains 
of compositions and sums of unbounded operators, we have that
\[ \dm(\Box)=\left\{f\in L^2_*(\Omega)\mid f\in\dm(\dbar)\cap\dm(\dbar^*) , \dbar f \in \dm(\dbar^*), \text{ and } \dbar^*f \in \dm(\dbar)\right\}. \]

A very important special case is when $\Omega$ is realized as a relatively compact and smoothly bounded domain
in a larger hermitian manifold $\mathcal{M}$, and given as $\Omega=\{z\in \mathcal{M}\mid \rho(z)<0\}$, where the gradient $\nabla\rho$ is normalized to unit length on
$\bdry \Omega$. In this case, if $f\in \mathcal{C}^2_*(\overline{\Omega})$ is a form 
smooth up to the boundary, the condition that $f\in\dm(\Box)$ is equivalent to $f$ satisfying on $\bdry \Omega$ the  {\em $\dbar$-Neumann boundary conditions}
\begin{equation}\label{eq-dbarneumann} \begin{cases}
f\rfloor \nabla \rho &=0, \text{ and }\\
\dbar f \rfloor\nabla \rho &=0,
\end{cases}
\end{equation}
where $\rfloor$ denotes the contraction of a form by a vector field.

\subsection{Differential forms on product manifolds}

We now generalize equations \eqref{eq-conten} and \eqref{eq-natten} to 
spaces of differential forms on manifolds.
Let $\h_1$ and $\h_2$ be vector spaces of differential forms on the manifolds 
$\Omega_1$ and $\Omega_2$ respectively. Let $\pi_j$ denote the projection from
the product $\Omega=\Omega_1\times\Omega_2$ to the factor $\Omega_j$.
 It is easy to see that the identification
\[ f\tensor g = \pi_1^*f \wedge \pi_2^* g\]
linearly extends to an injective map of $\h_1\tensor \h_2$ into  the space of differential forms
on $\Omega$. In particular, if we take $\h_j =L^2_*(\Omega_j)$, the Hilbert space of forms square integrable with 
respect to the hermitian metric (see \cite[Chapter~5]{cs} for detailed definitions) we get an injective map $L^2_*(\Omega_1)\tensor
 L^2_*(\Omega_2)\hookrightarrow L^2_*(\Omega)$ which can be extended 
by closure to obtain a natural identification
\[ L^2_*(\Omega_1)\csor L^2_*(\Omega_2)\cong L^2_*(\Omega).\]
Note that for $(p,q)$ forms, reading off degrees on each side, this construction gives a representation of  $L^2_{p,q}(\Omega)$
as an orthogonal direct sum of tensor products:
\begin{equation}\label{eq-lpq} L^2_{p,q}(\Omega)=\bigoplus_{\substack{p_1+p_2=p\\q_1+q_2=q}}L^2_{p_1,q_1}(\Omega_1)\csor L^2_{p_2,q_2}(\Omega_2).\end{equation}
It is clear how this construction extends to more than two factors.

\subsection{Proof of main theorem}\label{sec-mainthmproof}
 We introduce some more notation. Denote by $\Box^1$ and $\Box^2$ the 
complex Laplacians on $\Omega_1$ and $\Omega_2$ respectively (this is unambiguous since we never
consider powers of the complex Laplacian.) Then for $j=1,2$, the operator $\Box^j$ is
a densely defined, selfadjoint, nonnegative operator on $L^2_*(\Omega_j)$, and we we will denote 
its restriction to $L^2_{p,q}(\Omega_j)$ by $\Box_{p,q}^j$. Let $\Omega=\Omega_1\times \Omega_2$ be 
the product hermitian manifold (with product metric.) Let $D$ be the operator on
$L^2_*(\Omega)$ with domain $\dm(\Box^1)\tensor\dm(\Box^2)$, defined as
\begin{equation}\label{eq-D} D=\Box^1\tensor I_2 + I_1\tensor \Box^2,\end{equation}
where $I_j$ is the identity operator on $L^2_*(\Omega_j)$. Note that $D$ is densely defined on 
$L^2_*(\Omega)$, but it is not clear a priori whether $D$ is closable or not. 

Denote by $\Box$ the complex Laplacian on the product $\Omega$. We claim that $D$ is closable,
and its closure is $\Box$. The first step of the proof is the following:
\begin{lem}\label{lem-dbox}
$\dm(D)\subset\dm(\Box)$, and the restriction of $\Box$ to $\dm(D)=\dm(\Box^1)\tensor\dm(\Box^2)$ coincides with $D$.
\end{lem} 
\begin{proof} Denote by  $\dbar_j$ the $L^2$ Cauchy-Riemann operator
on $\Omega_j$ and by $\dbar^*_j$ its Hilbert adjoint, and let $\dbar, \dbar^*$ denote the corresponding objects
on the product $\Omega=\Omega_1\times\Omega_2$. It is easy to see from the definitions of the domains of these operators
that $\dm(\dbar)\supset \dm(\dbar_1)\tensor\dm(\dbar_2)$ and $\dm(\dbar^*)\supset \dm(\dbar_1^*)\tensor\dm(\dbar_2^*)$.
Further, on $\dm(\dbar_1)\tensor\dm(\dbar_2)$ we have the Leibniz formula
\begin{equation}\label{eq-leib} \dbar = \dbar_1\tensor I_1 + \sigma_1 \tensor \dbar_2,\end{equation}
where $\sigma_1$ is the operator on $L^2_*(\Omega_1)$, which when restricted to $L^2_{p,q}(\Omega_1)$ is
multiplication by $(-1)^{p+q}$. Note that $\sigma_1^2=I_1$ and for any operator $S$ of degree $d$ 
on $L^2_*(\Omega_1)$
(i.e., $\deg (Sf)-\deg (f)=d$ for every homogeneous form $f$ in $\dm(S)$) we have $\sigma_1 S = (-1)^d S \sigma_1$.

 Similarly for $\dbar^*$ we have on $\dm(\dbar_1^*)\tensor\dm(\dbar_2^*)$ that
\begin{equation}\label{eq-leibstar} \dbar^* = \dbar_1^*\tensor I_1 + \sigma_1 \tensor \dbar_2^*.\end{equation}

Now let $f_j\in \dm(\Box^j)$ and set $f= f_1\tensor f_2$. We verify that $f\in \dm(\Box)$. 
Indeed, since $f_j\in\dm(\dbar_j)\cap \dm(\dbar_j^*)$, it follows that 
\[ f\in \left(\dm(\dbar_1)\tensor\dm(\dbar_2)\right)\cap \left(\dm(\dbar_1^*)\tensor\dm(\dbar_2^*)\right).\]
Now, using also the facts that $\dbar f_j\in\dm(\dbar^*)$ and \eqref{eq-leib}, it follows that $\dbar f\in\dm(\dbar^*)$.
Similarly, using $\dbar^*f_j\in\dm(\dbar)$ and \eqref{eq-leibstar} we obtain that $\dbar^*f\in \dm(\dbar)$. It follows
that $f\in \dm(\Box)$. So $\dm(D)\subset\dm(\Box)$.

Now we compute
\begin{align*}
\dbar\dbar^* f &=(\dbar_1\tensor I_2 + \sigma_1\tensor \dbar_2)(\dbar_1^*f_1\tensor f_2 +\sigma_1f_1\tensor \dbar_2^*f_2)\\
&= \dbar_1\dbar_1^*f_1\tensor f_2
+ \dbar_1\sigma_1f_1\tensor\dbar_2^*f_2
+\sigma_1\dbar_1^*f_1\tensor\dbar_2f_2+
f_1\tensor\dbar_2\dbar_2^*f_2
\end{align*}
and
\begin{align*}
\dbar^*\dbar f&=(\dbar_1^*\tensor I_2+\sigma_1\tensor \dbar_2^*)(\dbar_1f_1\tensor f_2+\sigma_1 f_1\tensor \dbar f_2)\\
&=\dbar_1^*\dbar_1f_1\tensor f_2+\sigma_1\dbar_1f_1 \tensor\dbar_2^*f_2 +\dbar_1^*\sigma_1f_1\tensor \dbar_2 f_2+ f_1\tensor\dbar_2^*\dbar_2 f_2.
\end{align*}
Combining we obtain,
\begin{align*}
\Box f &= \Box^1 f_1\tensor f_2 +(\dbar_1\sigma_1+\sigma_1\dbar_1)f_1\tensor \dbar_2^*f_2
+(\sigma_1\dbar_1^*+\dbar_1^*\sigma_1)f_1\tensor \dbar_2f_2
+f_1\tensor \Box^2f_2\\
&= \Box^1f_1\tensor f_2+ f_1\tensor \Box^2f_2\\
&= Df,
\end{align*}
where we have made use of the fact that $\dbar$ and $\dbar^*$ are of  degree $\pm 1$ respectively.
By linear extension, it follows that $\Box=D$ on $\dm(D)$.
\end{proof}
The proof of Theorem~\ref{thm-specsum} will also require the following simple observation: 
\begin{lem}\label{lem-EF} Let $E$ and $F$ be closed subsets of the set  of nonnegative reals. Then
$E+F$ is also a closed set.
\end{lem}
\begin{proof} Let $z$ a point in the closure $\overline{E+F}$, and let $z_\nu$ be a sequence of points in
$E+F$ converging to $z$. Writing $z_\nu=x_\nu+y_\nu$, where $x_\nu\in E$ and $y_\nu\in F$ we see that $x_\nu$, $y_\nu$ are
bounded sequences in the closed sets $E$ and $F$ respectively. By the Bolzano-Weierstrass theorem, after 
passing to a subsequence, we can assume that $x_\nu\rightarrow x\in E$ and $y_\nu\rightarrow y\in F$. It follows 
that $z=x+y\in E+F$.
\end{proof}
\begin{proof}[Proof of equation \eqref{eq-sigma}]We first consider the case when $N=2$.
Since $\Box^1$ and $\Box^2$ are selfadjoint, using the results
of $\S$\ref{sec-AB}, we see that $D$ is an essentially self adjoint operator on $L^2_*(\Omega)$. By Lemma~\ref{lem-dbox},
$\Box$ is an extension of $D$. But since $\Box$ is selfadjoint, this means that $\Box$ is the closure of $D$. Thanks 
again to the results of $\S$\ref{sec-AB}, it follows that
\[ \spc(\Box)=\overline{\spc(\Box^1)+\spc(\Box^2)}.\]
Now, using Lemma~\ref{lem-EF}, and the fact that $\Box^j$ is nonnegative, so its spectrum consists of nonnegative numbers, it follows that
\[ \sigma(\Omega) = \sigma(\Omega_1) +\sigma(\Omega_2).\]
The case of general $N>2$ now follows by  a straightforward induction argument.

\end{proof}
\begin{proof}[Proof of equation \eqref{eq-sigmapq}] First we assume that $N=2$.  Let $f_1\in L^2_{p_1,q_1}(\Omega_1)$,
and $f_2\in L^2_{p_2,q_2}(\Omega_2)$, so that $f_1\tensor f_2\in L^2_{p_1,q_1}(\Omega_1)\csor L^2_{p_2,q_2}(\Omega_2)$.
Note that $\Box^1$ (resp. $\Box^2$) maps $L^2_{p_1,q_1}(\Omega_1)$ ( resp. $L^2_{p_2,q_2}(\Omega_2)$) into itself.
Therefore, the formula \eqref{eq-D} defining $D$ shows that $D$ maps  $L^2_{p_1,q_1}(\Omega_1)\csor L^2_{p_2,q_2}(\Omega_2)$
into itself, and it follows that so does $\Box$, the closure of $D$. Therefore, the restriction of $\Box$ defines a selfadjoint operator on each space $ L^2_{p_1,q_1}(\Omega_1)\csor L^2_{p_2,q_2}(\Omega_2)$. Denote this restriction by
the admittedly barbarous notation $\Box_{\substack{p_1,q_1\\p_2,q_2}}$. Then by the results of $\S$\ref{sec-AB}, we have that $\Box_{\substack{p_1,q_1\\p_2,q_2}}$ is the unique self adjoint extension of the operator $\Box^1_{p_1,q_1}\tensor I_2+ I_1\tensor \Box^2_{p_2,q_2}$ (where now $I_j$ denotes the identity map on the Hilbert space $L^2_{p_j,q_j}(\Omega_j)$,) and we have for the spectra
\begin{align}
\spc\left(\Box_{\substack{p_1,q_1\\p_2,q_2}}\right)&= \overline{\spc(\Box^1_{p_1,q_1})+\spc(\Box^2_{p_2,q_2})}\nonumber\\
&=\overline{\sigma_{p_1,q_1}(\Omega_1)+\sigma_{p_2,q_2}(\Omega_2)}\nonumber\\
&= \sigma_{p_1,q_1}(\Omega_1)+\sigma_{p_2,q_2}(\Omega_2)\label{eq-boxspec},
\end{align}
where we have used Lemma~\ref{lem-EF}.

Now consider a Hilbert space $\h$ represented as an orthogonal direct sum
\[ \h = \bigoplus_{k=1}^n \h_k,\]
and let $A$ be an operator on $\h$ which maps each $\h_k$ to itself. Denoting by
$A_k$ the restriction of $A$ to $\h_k$ (interpreted as an operator on $\h_k$,)
there is a direct sum decomposition
$A= \bigoplus_{k=1}^n A_k$, Then we  have (see \cite[Theorem~2.23]{tes}):
\begin{equation}\label{eq-sumspec} \spc(A)= \bigcup_{k=1}^n \spc(A_k).\end{equation}

Now, by \eqref{eq-lpq}, the space $L^2_{p,q}(\Omega)$ is represented as a orthogonal direct sum 
of subspaces $L^2_{p_1,q_1}(\Omega_1)\csor L^2_{p_2,q_2}(\Omega_2)$ (with $p_1+p_2=p$ and $q_1+q_2=q$,)
and the operator $\Box_{p,q}$ maps each of these subspaces to itself. Therefore, using 
\eqref{eq-sumspec}, we have
\begin{align*}
\spc(\Box_{p,q})&= \bigcup_{\substack{p_1+p_2=p\\q_1+q_2=q}}\spc\left(\Box_{\substack{p_1,q_1\\p_2,q_2}}\right)\\
&=\bigcup_{\substack{p_1+p_2=p\\q_1+q_2=q}}\left(\spc(\Box^1_{p_1,q_1})+\spc(\Box^2_{p_2,q_2})\right)&\text{ by \eqref{eq-boxspec}}
\end{align*}
In the notation used in the statement of Theorem~\ref{thm-specsum}, this reads
\[ \sigma_{p,q}(\Omega) = \bigcup_{\substack{p_1+p_2=p\\q_1+q_2=q}}\left(\sigma_{p_1,q_1}(\Omega_1)+ \sigma_{p_2.q_2}(\Omega_2)\right),\]
which proves  \eqref{eq-sigmapq} in the case $N=2$.

Now we extend this result by induction to general $N>2$. The induction will require the use of
the following formula regarding Minkowski sums --
\begin{equation}\label{eq-mindist}
E + \bigcup_{i=1}^n F_i = \bigcup_{i=1}^n\left(E+F_i\right),
\end{equation}
which follows directly from the definition. We assume that the result has been established for 
$N-1$ factors, and consider $N$ factors $\Omega_1,\dots, \Omega_N$. We set $\Omega_j' =\Omega_j$ for $1\leq j\leq N-2$,
and let $\Omega_{N-1}' =\Omega_{N-1}\times \Omega_N$. Then, if $\Omega= \Omega_1\times\dots\times\Omega_N=\Omega_1'\times\dots\times\Omega_{N-1}'$, then we have by the induction hypothesis
\begin{align*}
 \sigma_{{\boldsymbol{p}},{\boldsymbol{q}}}(\Omega)& =\bigcup_{\substack{\sum_{j=1}^{N-1} p_j={\boldsymbol p}\\\sum_{j=1}^{N-1} q_j={\boldsymbol q}}}
\left(\sigma_{p_1,q_1}(\Omega_1')+\dots+\sigma_{p_{N-1},q_{N-1}}(\Omega_{N-1}')\right)\\
&=  \bigcup_{\substack{\sum_{j=1}^{N-1} p_j={\boldsymbol p}\\\sum_{j=1}^{N-1} q_j={\boldsymbol q}}}
\left(
	\sigma_{p_1,q_1}(\Omega_1)+\dots+\sigma_{p_{N-2},q_{N-2}}(\Omega_{N-2}) + 
       \left(\bigcup_{\substack{P_1+P_2=p_{N-1}\\Q_1+Q_2=q_{N-1}}}\left(\sigma_{P_1,Q_1}(\Omega_{N-1})+\sigma_{P_2,Q_2}(\Omega_{N})\right)\right)\right)\\
&\text{(using the result for $N=2$)}\\
&=  \bigcup_{\substack{\sum_{j=1}^{N-1} p_j={\boldsymbol p}\\\sum_{j=1}^{N-1} q_j={\boldsymbol q}}}
\left(\bigcup_{\substack{P_1+P_2=p_{N-1}\\Q_1+Q_2=q_{N-1}}}\left(
\sigma_{p_1,q_1}(\Omega_1)+\dots+\sigma_{p_{N-2},q_{N-2}}(\Omega_{N-2}) +\sigma_{P_1,Q_1}(\Omega_{N-1})+\sigma_{P_2,Q_2}(\Omega_{N})\right)\right)\\
&\text{(using \eqref{eq-mindist})}\\
&=\bigcup_{\substack{\sum_{j=1}^N p_j={\boldsymbol p}\\\sum_{j=1}^N q_j={\boldsymbol q}}}
\left(\sigma_{p_1,q_1}(\Omega_1)+\dots+\sigma_{p_N,q_N}(\Omega_N)\right)\\
&\text{(renaming the indices.)}
\end{align*}
This completes the proof.
\end{proof}

\section{Proof of Theorem~\ref{thm-zucker}}
We first establish a couple of lemmas :
\begin{lem}\label{lem-tfae}
Let $A$ be a nonnegative selfadjoint operator on a Hilbert space $\h$. Then the following are
equivalent:
\begin{enumerate}
\item The range of $A$ is closed.
\item There is a $C>0$ such that for each $x\in \dm(A)\cap \ker(A)^\perp$
\begin{equation}\label{eq-closedran}
\norm{Ax}\geq C\norm{x}.
\end{equation}
\item There is a $c>0$ such that the intersection of the open interval
$(0,c)$ with the set $\spc(A)$ is empty.
\end{enumerate}
\end{lem}
\begin{proof}
The equivalence of (1) and (2) is a standard fact in functional analysis (see \cite[Lemma~4.1.1]{cs}.)
We show that (2) and (3) are equivalent. Using the spectral theorem as stated in $\S$\ref{sec-spthm}, we can assume that $\h$ is the space $L^2(X,\mu)$ for some measure space $(X,\mu)$, and $Af= hf$ for a nonnegative function $h$ on $X$.

First assume that (3) is true. Recall from $\S$\ref{sec-multop} that the spectrum of $A$ 
coincides with the essential range of the function $f$. This means that on the complement of the set $\{h=0\}\subset X$,
the function $h$ satisfies $h\geq c$ a.e.. Then, we have for any $f\in\dm(A)\subset L^2(X,\mu)$,
\begin{align*}
\norm{Af}^2&=\int \abs{hf}^2 d\mu\\
&=\int_{\{h=0\}}\abs{hf}^2 d\mu +\int_{\{h>0\}} \abs{hf}^2d\mu\\
&\geq0+ c^2\int \abs{f}^2 d\mu\\
&= c^2 \norm{f}^2,
\end{align*}
so that (3) implies (2). Now assume that (3) is violated. For a positive integer $\nu$, let
\[ E_\nu=\left\{x\in X\mid \frac{1}{2^{\nu+1}}\leq h(x) <\frac{1}{2^{\nu}}\right\}\]
and let $\mu_\nu= \mu(E_\nu)$. Since (3) is not true, it is possible to find a sequence 
of integers $\nu_k\uparrow\infty$ such that $\mu_{\nu_k}>0$. For each $k$ define
\[ f_k = \begin{cases} \displaystyle\frac{1}{\sqrt{\mu_{\nu_k}}} &\text{ on $E_{\nu_k}$}\\
			0 & \text{elsewhere.}\\
	\end{cases}
\]
We claim that $f_k$ is orthogonal to $\ker(A)$. Indeed, $g\in \ker(A)$ if and only if $g$ has support in the set 
$\{h=0\}$ (cf. equation \eqref{eq-eigen}.) Since the support of $f_k$ is by construction disjoint from that of $g$, it follows that
$\int f_k\overline{g} d\mu=0$.
Also, 
\begin{align*}\norm{f_k}^2&= \int_{E_{\nu_k}}\frac{1}{\mu_{\nu_k}}d\mu\\
&=1,
\end{align*}
but on the other hand we have
\begin{align*}
\norm{Af_k}^2 &=\int \abs{hf_k}^2d\mu\\
&=\int_{E_{\nu_k}}\frac{h^2}{\mu_{\nu_k}}d\mu\\
&\leq \frac{1}{4^{\nu_k}}.
\end{align*}
It follows that $f_k\in\dm(A)\cap \ker(A)^\perp$, but no  constant such as $C$ in \eqref{eq-closedran} exists. This completes the proof.
\end{proof}
\begin{lem} \label{lem-nonnegative} Let $A_1$, $A_2$ be nonnegative selfadjoint operators on Hilbert spaces $\h_1, \h_2$,
 and let $A$ 
be the closure of  $A_1\tensor I_2+I_1\tensor A_2$ as an operator on $\h_1\csor\h_2$. (Recall that $A$  was shown 
in $\S$\ref{sec-AB} to be selfadjoint.) Then we have 
\begin{equation}\label{eq-kernel} \ker(A) =\ker(A_1)\csor\ker(A_2).\end{equation}
\end{lem}
\begin{proof}
We use the representation of $A_1,A_2,A$ by multiplication operators developed in $\S$\ref{sec-AB}.
After using the unitary isomorphism $U$, proving \eqref{eq-kernel} is reduced to proving that 
$\ker(T_h)=\ker(T_{h_1})\csor \ker(T_{h_2})$. 
The functions $h_1$ and $h_2$ on $X_1$ and $X_2$ respectively
which represent $A_1$ and $A_2$ by multiplication are now nonnegative a.e.,
and the subset $h^{-1}(0)$ of $X_1\times X_2$ is identical to $h_1^{-1}(0)\times h_2^{-1}(0)$.
Indeed, since $h_1\geq 0$ and $h_2\geq 0$ a.e., the only way $h(x_1,x_2)= h_1(x_1)+h_2(x_2)$ can vanish a.e.
is by the vanishing a.e. of both $h_1$ and $h_2$. Using the representation \eqref{eq-eigen} of the kernel,
we obtain
\begin{align*}
\ker(T_h) &= L^2(h^{-1}(0),\mu_1\tensor\mu_2)\\
&=L^2(h_1^{-1}(0)\times h_2^{-1}(0), \mu_1\tensor \mu_2)\\
&= L^2(h_1^{-1}(0),\mu_1)\csor L^2(h_2^{-1}(0),\mu_2)\\
&=\ker(T_{h_1})\csor\ker(T_{h_2})
\end{align*}
\end{proof}
\begin{proof}[Proof of Theorem~\ref{thm-zucker}] It is sufficient to prove the result for $N=2$,
the general case following by a simple induction argument. We recall the following standard fact from Hodge theory:
on a hermitian manifold the following are equivalent:
(a) $\Box$ has closed range, 
(b) $\dbar$ has closed range,
(c) $\dbar^*$ has closed range,
(d) (Kodaira; see \cite[p.~165]{dRh}, or \cite[Lemma~2.2]{chsh}) every $L^2$-Dolbeault cohomology class
has a unique harmonic representative: more precisely,  the inclusion $\ker(\Box)\subset \ker(\dbar)$ 
induces as isomorphism on the cohomology level.

We use the notation used in the proof 
of Lemma~\ref{lem-dbox}. Therefore, assume that the operators  $\dbar_1$ and $\dbar_2$ have closed range
in $L^2_*(\Omega_1)$ and $L^2_*(\Omega_2)$ respectively. Therefore for $j=1,2$, the operator and $\Box^j$ also
has closed range in $L^2_*(\Omega_j)$. It follows from Lemma~\ref{lem-tfae} that there exist $c_j>0$ such that
$\sigma(\Omega_j)\cap(0,c_j)=\emptyset$. But $\sigma(\Omega)=\sigma(\Omega_1)+\sigma(\Omega_2)$ by 
Theorem~\ref{thm-specsum}, so if $c=\min(c_1,c_2)$, then clearly $\sigma(\Omega)\cap(0,c)=\emptyset$. It follows
now from Lemma~\ref{lem-tfae} that $\Box$ has closed range in $L^2_*(\Omega)$. 
We conclude that $\dbar$ has closed range in $L^2_*(\Omega)$.

Now, $\Box^j$ is a nonnegative operator, $j=1,2$, and $\Box$ is the unique selfadjoint extension of
$\Box^1\tensor I_2+I_1\tensor \Box^2$, so by Lemma~\ref{lem-nonnegative},  we have
\[ \ker(\Box) =\ker(\Box^1)\csor\ker(\Box^2).\]
Since, by part (d) of the result quoted in the first paragraph of this proof,
$\ker(\Box)\cong H^*_{L^2}(\Omega)$ and $\ker(\Box^j)\cong H^*_{L^2}(\Omega_j)$ for $j=1,2$,
the K\"{u}nneth formula \eqref{eq-kun} now follows in the case $N=2$. 
\end{proof}

\section{Point spectra and eigenvectors}
In this section we consider certain simple special cases of Theorem~\ref{thm-specsum}. All these could have 
been deduced directly from the representation of $\Box$ constructed in $\S$\ref{sec-mainthmproof}. However,
we give direct elementary arguments wherever possible in view of the importance of the special cases considered.
\subsection{Expansions in eigenforms}
We use the notation of Theorem~\ref{thm-specsum} and  $\S$\ref{sec-mainthmproof}. 
\begin{prop}(a) For $j=1,\dots, N$,  let $\alpha_j\in\sigma(\Omega_j)$ be an eigenvalue, and let
$E_{\alpha_j}^j\subset\dm(\Omega_j)\subset L^2_*(\Omega_j)$ be the corresponding eigenspace. Then $\sum_{j=1}^N\alpha_j\in \sigma(\Omega)$ is
an eigenvalue and the corresponding eigenspace is $E_{\alpha_1}^1\csor\dots E_{\alpha_N}^N$.

(b) If each  $\sigma(\Omega_j)$ consists only of only eigenvalues, so does $\sigma(\Omega)$.
\end{prop}
\begin{proof}
By a simple induction argument, it suffices to consider the case $N=2$ for both (a) and (b). For 
part (a), we let $f_j\in E_{\alpha_j}^j\subset\dm(\Box_j)$ so that $f_1\tensor f_2 \in \dm(\Box)$ by Lemma~\ref{lem-dbox}.
A computation using \eqref{eq-D} and Lemma~\ref{lem-dbox}
now shows that $\alpha_1+\alpha_2$ is an eigenvalue of $\Box$ with eigenvector
$f_1\tensor f_2$. Part (a) follows, since the algebraic tensor product is dense in the 
Hilbert tensor product (cf. $\S$\ref{sec-tensor}.)

For part (b), continuing to assume $N=2$, we note that the hypothesis implies that for $j=1, 2$,
\[ L^2_*(\Omega_j) = \bigoplus_{\lambda\in \sigma(\Omega_j)} E^j_\lambda,\]
where $E^j_\lambda$ denotes the eigenspace of $\Box_j$ corresponding to the eigenvalue $\lambda$.
Therefore, we have
\[ L^2_*(\Omega) = \bigoplus_{\substack{\lambda\in \sigma(\Omega_1)\\\mu\in\sigma(\Omega_2)}} E^1_\lambda\csor E^2_\mu.\]
Therefore, the span of the the eigenspaces corresponding to the points of $\sigma(\Omega_1)+\sigma(\Omega_2)$
is dense in the Hilbert space $L^2_*(\Omega)$. It follows that the full spectral decomposition of $\Box$ on $\Omega$ is given by projection on the eigenspaces
corresponding to  $\sigma(\Omega_1)+\sigma(\Omega_2)$. Part (b) now follows.
\end{proof}

We now consider  $\Omega_j$,  $j=1,\dots, N$ such that each $\sigma(\Omega_j)$ consists of
eigenvalues only. This happens in many important cases, for example, when each $\Omega_j$ is 
a smoothly bounded pseudoconvex domain of  finite type  in some $\cx^{n_j}$.
For  $\lambda\in \sigma(\Omega_j)$ denote by
$\pi^j_\lambda$ the orthogonal projection from $L^2_*(\Omega_j)$ to the subspace $E^j_\lambda$. Then
we have the following representation 
\[ \Box^j =\sum_{ \lambda\in \sigma(\Omega_j)}\lambda \pi^j_\lambda,\]
where the series on the right converges in the strong operator topology, i.e. for every $f\in \dm(\Box^j)$,
$\sum_{ \lambda\in \sigma(\Omega_j)}\lambda \pi^j_\lambda f$ converges to $\Box^jf$ in the norm topology of $L^2_*(\Omega_j)$. Our computations show that:
\begin{cor}\label{cor-box} On $\Omega$:

\[ \Box =\sum_{\substack{\lambda_j\in\sigma(\Omega_j)\\j=1,\dots,N}} (\lambda_1+\dots+\lambda_N)\:\pi_{\lambda_1}^1\csor\dots\csor \pi_{\lambda_N}^N\]

\end{cor}
Here $\pi_{\lambda_1}^1\csor\dots\csor \pi_{\lambda_N}^N$  is the projection operator on $L^2_*(\Omega)$ obtained as the closure of the bounded operator $\pi_{\lambda_1}^1\tensor\dots \tensor\pi_{\lambda_N}^N$
on the algebraic tensor product subspace $L^2_*(\Omega_1)\tensor\dots\tensor L^2_*(\Omega_N)$.
This is clearly a projection onto a subspace of the eigenspace of $\Box$  corresponding to the eigenvalue 
$\lambda=\sum_j\lambda_j$ of $\Box$ (it is not necessarily the {\em  full} projection corresponding
to $\lambda$, since there may be more than one way of representing $\lambda$ as a sum of eigenvalues of  the complex Laplacian on the factor domains.)

Denote by $\pi^{j,(p,q)}_\lambda$ the projection from $L^2_{p,q}(\Omega_j)$ onto the eigenspace of $\Box_{p,q}$ corresponding to the eigenvalue $\lambda\in \sigma_{p,q}(\Omega_j)$. 
An argument similar to the one above shows also that
\begin{cor}\label{cor-boxpq} We also have on $\Omega$:
\[ \Box_{\boldsymbol{p,q}}=\bigoplus_{\substack{\substack{\sum_{j=1}^N p_j={\boldsymbol p}\\\sum_{j=1}^N q_j={\boldsymbol q}}}}\:
\sum_{\lambda_j\in\sigma_{p_j,q_j}(\Omega_j)}(\lambda_1+\dots+\lambda_N)\:\pi^{1,(p_1,q_1)}_{\lambda_1}\csor\dots\csor\pi^{N,(p_N,q_N)}_{\lambda_N}\]
\end{cor}
\subsection{Example: Polydomains} We now consider the special case in which each $\Omega_j$ is a bounded
domain in the complex plane $\cx$ with smooth boundary, so that  the pseudoconvex domain $\Omega=\Omega_1\times\dots\times\Omega_N\subset\cx^N$ is a so called {\em polydomain}. We consider the spectrum of $\Box_{0,q}$ on $\Omega$.
For convenience we will write
$\Box_q$ for $\Box_{0,q}$, and  $\sigma_q(\Omega_j)$ for $\sigma_{0,q}(\Omega_j)$.

Note that the spectrum $\sigma(\Omega_j)$ of $\Box^j$ consists of eigenvalues only.
Indeed $\Box^j_1$ is the same as the usual Laplacian with Dirichlet boundary conditions,
which is well-known to have a compact inverse, the  Green operator $\mathcal{G}$.
We can write the eigenvalues in $\sigma_1(\Omega_j)$ as an increasing sequence
\[  0<\mu_1^j\leq \mu_2^j \leq \dots\]
where we repeat each eigenvalue according to its (finite) multiplicity, and let
$Y_k^j(z_j)d\overline{z_j}$ denote an eigenform of $\Box^j_1$ corresponding to the eigenvalue $\mu_k^j$,
where $z_j$ denotes the natural coordinate on $\Omega_j$, and the eigenforms are so chosen for
each eigenvalue with multiplicity that
the collection $\{Y_k^j\}_{k\in\mathbb{N}}$ is a complete orthogonal set in $L^2(\Omega_j)$

Now $u\in\Box^j_0\cap\mathcal{C}^1(\overline{\Omega_j})$
means that $\dbar u=0$ on $\bdry \Omega_j$, by \eqref{eq-dbarneumann}. For smooth $f$ therefore,
the equation $\Box^j_0 u=f$ takes the form
 \[ \begin{cases} \Delta u=-4 f  &\text{on $\Omega_j$}\\
		\displaystyle{\frac{\partial u}{\partial \overline{z}}}=0 & {\text{on $\bdry \Omega_j$}}\\
\end{cases}
\]
If $v=\frac{\partial u}{\partial \overline{z}}$, this can be rewritten as the Dirichlet problem 
\[ \begin{cases} \Delta v = -4f_{\overline{z}}&\text{on $\Omega_j$} \\
		v=0 & \text{on $\bdry \Omega_j$},\\
\end{cases}
\]	
where $\Delta$ is the usual Laplacian on $\rl^2$. It is now easy to see that 
\[ u = \frac{1}{\pi z}\ast \left(\mathcal{G} (-4f_{\overline{z}})\right)^0+ h\]
where $h$ is a $L^2$ holomorphic function, $\mathcal{G}$ is the (compact) solution operator 
of the Dirichlet Laplacian, $g\mapsto g^0$ is the extension-by-zero of a function on $\Omega_j$ to $\cx$,
and $\frac{1}{\pi z}$ is the fundamental solution of the $\dbar$-equation on $\cx$. It easily
follows that the inverse modulo the kernel of $\Box^j_0$ is compact, 
i.e. the restriction of $\Box^j_0$ to $\ker(\Box^j_0)^\perp$ has compact inverse. 
$\sigma_0(\Omega_j)$, the spectrum $\Box^j_0$ consists of eigenvalues only, which
can be written in ascending order as
\[ 0=\lambda_0^j<\lambda_1^j\leq \lambda_2^j\leq \dots,\]
where the positive eigenvalues are of finite multiplicity and they are repeated according to their  multiplicity. As noted above,
the eigenspace corresponding to the eigenvalue $\lambda_0^j=0$ is the Bergman space of $L^2$-holomorphic functions
on $\Omega_j$, and let $\{H_k^j\}_{k\in \mathbb{N}}$ be a complete orthogonal set in the Bergman space $L^2(\Omega_j)\cap \mathcal{O}(\Omega_j)$.
For $k\geq 1$, let $Z_k^j$ be an eigenfunction  of $\Box^j_0$ corresponding to  eigenvalue $\lambda_k^j$. We can again
assume that these have been chosen such that the family $\{Z_k^j\}_{k\in\mathbb{N}}$ is a complete orthogonal set
in $\ker(\Box^j_0)^\perp$.

For a subset $J$ of $\{1,\dots,n\}$ of cardinality $q$, where $J=\{j_1,\dots,j_q\}$ with 
$j_1<\dots<j_q$, we write $d\overline{z}^J=d\overline{z}_{j_1}\wedge\dots\wedge d\overline{z}_{j_q}$, with the 
understanding that $d\overline{z}^\emptyset =1$. We also use the standard convention that a sum over an 
empty set is 0.  With these notational preliminaries, Corollary~\ref{cor-boxpq} gives
rise to the following description of the eigenstructure of the operator $\Box_q$ on $\Omega$:
\begin{prop}\label{prop-polydomain}
Let $J$ be a subset of $\{1,\dots, n\}$ of cardinality $q$, and let $\mathbf{k}=(k_1,\dots,k_n)\in \mathbb{N}_+^n$ be 
an $n$-tuple of positive integers. Then

\[ \mu(J,\mathbf{k})= \sum_{j\in J}\mu^j_{k_j}\]
is an eigenvalue of $\Box_q$, and
\[ W(J,\mathbf{k})= \left(\prod_{j\in J}Y_{k_j}^j(z_j)\prod_{j\not\in J}H^j_{k_j}(z_j)\right)d\overline{z}^J\]
is an eigenform corresponding to this eigenvalue. Further, if $q<n$, then
\[\lambda(J,\mathbf{k})= \sum_{j\in J}\mu^j_{k_j}+\sum_{j\not\in J}\lambda^j_{k_j}\]
is also an eigenvalue of $\Box_q$, with eigenform
\[ V(J,\mathbf{k})=\left(\prod_{j\in J}Y_{k_j}^j(z_j)\prod_{j\not\in J}Z^j_{k_j}(z_j)\right) d\overline{z}^J.\]

Moreover, this is the complete list of eigenvalues and eigenforms of $\Box_q$ as $J$ ranges over all subsets 
of $\{1,\dots,n\}$ of size $q$ and $\mathbf{k}$ ranges over $ \mathbb{N}_+^n$, and gives the full spectral decomposition
of $\Box_q$.

\end{prop}
If $q<n$, the eigenvalue $\mu(J,\mathbf{k})$ has infinite multiplicity, since there are  
infinity many $\mathbf{k}$ corresponding to the same eigenvalue, and for distinct $\mathbf{k}$ we 
have distinct eigenforms $W(J,\mathbf{k})$. If $q=n$ on the other hand,  all the eigenvalues $\mu(J,\mathbf{k})$
are of finite multiplicity, as one would expect 
from the Dirichlet problem in a bounded domain.  Since $\Box_q$ has eigenvalues of infinite multiplicity
for $q<n$, it immediately follows that for $0<q<n$, the inverse of $\Box_q$, the $\dbar$-Neumann operator
$N_q$ is non-compact.

The special case of Proposition~\ref{prop-polydomain} when $\Omega_j=\{z\in \cx\mid \abs{z}<a_j\}$ for some $a_j>0$ 
(so that $\Omega$ is a polydisc)
was obtained in the paper \cite{fu}. In this case, the functions
$Y^j_k$ and $Z^j_k$ have explicit representations in terms of Bessel functions.

%********************************************************************


\begin{thebibliography}{12}
\bibitem{jb}Bertrams, Julia:
 Randregularit\"{a}t von L\"{o}sungen der $\overline\partial$ -Gleichung auf dem Polyzylinder und zweidimensionalen analytischen Polyedern.  Dissertation, Rheinische Friedrich-Wilhelms-Universit\"{a}t Bonn, Bonn, 1986.

\bibitem{brun} Br\"{u}ning, J; Lesch, M.:
 Hilbert complexes.
 {\em J. Funct. Anal.} {\bf 108}  (1992),  no. 1, 88--132.
\bibitem{chsh} Chakrabarti, Debraj and Shaw, Mei-Chi: The Cauchy-Riemann equations 
on a product domain. Preprint. Available at \texttt{arXiv:0911.0103}
\bibitem{chgr2}  Cheeger, Jeff: On the Hodge theory of Riemannian pseudomanifolds. {\em Proc. Sympos. Pure Math.}
{\bf 36}, Amer. Math. Soc, Providence, RI, (1980), 91-145.
\bibitem{chgr1} Cheeger, Jeff: Spectral geometry of singular Riemannian spaces.
{\em J. Differential Geom.} {\bf 18}, no. 4 (1983), 575-657. 
\bibitem{cs}
Chen, So-Chin and Shaw, Mei-Chi:
{\em  Partial differential equations in several complex variables.}
AMS/IP Studies in Advanced Mathematics, {\em 19.}
American Mathematical Society, Providence, RI;
International Press, Boston, MA, 2001.
\bibitem{dRh}de Rham, Georges.
{\em Vari\'{e}t\'{e}s diff\'{e}rentiables. Formes, courants, formes harmoniques.}  Hermann et Cie, Paris, 1955.
\bibitem{eh1} Ehsani, Dariush:
 Solution of the $\overline\partial$-Neumann problem on a non-smooth domain.
 {\em Indiana Univ. Math. J.}{\bf  52 } (2003),  no. 3, 629--666.
\bibitem{eh2} Ehsani, Dariush:
 Solution of the $\overline\partial$-Neumann problem on a bi-disc.
{\em  Math. Res. Lett.}{\bf   10 } (2003),  no. 4, 523--533.
\bibitem{eh3}Ehsani, Dariush:
 The $\overline\partial$-Neumann problem on product domains in $\cx^n$.
{\em  Math. Ann.}  {\bf 337}  (2007)  797--816.


\bibitem{fk} 
Folland, Gerald   B. and  Kohn, Joseph  J.:   \newblock
{\em The Neumann problem for the Cauchy-Riemann complex.}
Annals of Mathematics Studies, No. {\bf 75.}
Princeton University Press, Princeton, N.J.;
University of Tokyo Press, Tokyo, 1972.
\bibitem{fu} Fu, Siqi:
Spectrum of the $\overline\partial$-Neumann Laplacian on polydiscs. 
{\em Proc. Amer. Math. Soc.} {\bf 135} (2007), no. 3, 725--730.
\bibitem{hal}Halmos, P. R: What does the spectral theorem say?
{\em Amer. Math. Monthly} {\bf 70} (1963) 241--247. 
\bibitem{rs}Reed, Michael and Simon, Barry:
{\em Methods of modern mathematical physics, Vol. I:
Functional analysis.} Second edition. Academic Press, Inc., New York, 1980.
\bibitem{mcs1}Shaw, Mei-Chi: Global solvability and regularity for $\dbar$ on an annulus 
between two weakly pseudoconvex domains. {\em Trans.
 Amer. Math. Soc.}{\bf 291} (1985), 255–267.
\bibitem{mcs2}Shaw, Mei-Chi: The closed range property for ∂ on domains with pseudoconcave boundary. 
{\em Proceedings of the Complex Analysis conference, Fribourg, Switzerland} (2008).
\bibitem{tes} Teschl, Gerald: {\em Mathematical methods in quantum mechanics.
With applications to Schrödinger operators.} Graduate Studies in Mathematics, {\bf 99.} 
American Mathematical Society, Providence, RI, 2009. 
\bibitem{zuck} Zucker, Steven:  $L_{2}$ cohomology of warped products
and arithmetic groups.{\em   Invent. Math.}{\bf  70 } (1982),  no. 2, 169--218.
\end{thebibliography}
\end{document}